\documentclass[11pt,a4paper,reqno]{amsart}
\usepackage{float}
\usepackage{amsfonts,amsmath,graphics,epsfig,latexsym,times}
\usepackage[all]{xy}
\newif\ifpix \pixtrue
\usepackage[english]{babel}

\numberwithin{equation}{section}

\usepackage{check}

\newcommand{\chM}{\check M}
\renewcommand{\chM}{\goodMcheck}

\dedicatory{For Professor Paul Gauduchon on his 60th birthday}

\calclayout

\numberwithin{equation}{section}

\def\TT{{\mathbb{T}}} \def\PP{{\mathbb{P}}} \def\HH{{\mathbb{H}}}
\def\ZZ{{\mathbb{Z}}} \def\QQ{{\mathbb{Q}}} \def\CC{{\mathbb{C}}}
\def\RR{{\mathbb{R}}} \def\EE{{\mathbb{E}}} 
\def\cO{{\mathcal{O}}}
\def\cX{{\mathcal{X}}}

\def\cU{{\mathcal{U}}}

\renewcommand{\epsilon}{\varepsilon}

\newcommand{\SU}{\mathrm{SU}}
\newcommand{\zb}{\overline{z}}

\newcommand{\SO}{\mathrm{SO}}

\newcommand{\CP}{\mathbb{CP}}

\newcommand{\del}{\partial}
\newcommand{\delb}{\bar\partial}

\newcommand{\ovS}{\overline{\Sigma}}

\newcommand{\ovpi}{\overline{\pi}}

\newcommand{\ovM}{\overline{M}}

\newcommand{\hM}{\widehat{M}}

\newcommand{\hS}{\widehat{\Sigma}}

\newcommand{\rd}{{\mathrm d}}

\newtheorem{lemma}[subsubsection]{Lemma}
\newtheorem{prop}[subsubsection]{Proposition}

\newtheorem{theo}[subsection]{Theorem}
\newtheorem{theointro}{Theorem}

\theoremstyle{definition}

\newtheorem{question}[subsection]{Question}
\newtheorem*{rmk*}{Remark}
\newtheorem*{rmks*}{Remarks}
\newtheorem{rmks}[subsection]{Remarks}
\renewenvironment{proof}{\paragraph{\emph{Proof}}} {~Q.E.D.\medskip}
\newenvironment{remark}{\subsubsection{\bf Remark}--- \normalfont} { }
\newenvironment{remark*}{\begin{rmk*} --- \normalfont} { \end{rmk*} }
\newenvironment{remarks*}{\begin{rmks*} \begin{enumerate} \normalfont}
{\end{enumerate} \end{rmks*} } 
 
\title[K\"ahler surfaces of constant scalar curvature]{Construction of  K\"ahler surfaces with constant scalar
  curvature}
\date{December 2004}
\author{Yann Rollin}
\thanks{First author partly supported by NSF grant \# DMS-0305130}
\address{Yann Rollin, MIT, 77 Massachusetts Avenue,
Cambridge MA 02139, USA}
\email{rollin@math.mit.edu}
\author{Michael Singer}
\thanks{Second author partly supported by a Leverhulme research
fellowship and an EPSRC small grant}  
\address{Michael Singer, School of Mathematics, King's buildings,
  Edinburgh Scotland} 
\email{michael@maths.ed.ac.uk}
\thanks{}
\begin{document}
{\Huge \sc \bf\maketitle}
\begin{abstract}
We present new constructions of K\"ahler metrics with constant scalar
curvature on complex surfaces, in particular on certain del Pezzo
surfaces. Some higher dimensional examples are provided as well.
\end{abstract}
\section{Introduction}
\label{secintro}
The aim of this note is to present a  new construction of K\"ahler metrics of
constant scalar curvature (CSC) on complex surfaces. In order to introduce our results, let us introduce the terms ``positive CSC'' to mean ``constant positive scalar curvature'', ``zero CSC'' for ``(constant) zero scalar curvature'' and ``negative CSC'' for ``constant negative scalar curvature''.

Our construction
gives rise to many families of examples, but in this introduction we
shall focus on 
$$
X_k := \mbox{$k$-fold blow-up of }\CP^1\times \CP^1.
$$
We note that if $k\geq 1$ then $X_k$ can also be viewed as a
 $k+1$-fold blow-up of $\CP^2$.  Of course, the above description of $X_k$ does not fix its complex structure: this will depend on the location of the centres of the blow-ups.  Our first result gives positive CSC K\"ahler metrics in a family of K\"ahler classes on $X_k$, for $k=6, 7, 8$, and for certain choices of complex structure.  We note that if $k\leq 7$ then $X_k$ is Fano and the work of Tian \cite{T} and others gives positive K\"ahler--Einstein metrics on $X_k$.  Our result is new in that it produces CSC metrics on $X_8$ as
well as CSC metrics on $X_6$ and $X_7$ in K\"ahler classes that are
``arbitrarily far'' from $c_1(X)$:
\begin{theointro}
\label{cp9}
For $k=6,7,8$, there exists a $k$-point blow-up $X$ of
$\CP^1\times\CP^1$ with no non-trivial holomorphic vector field and the
following properties. Let $F=\{x\}\times \CP^1$ be a generic rational curve of
$\CP^1\times\CP^1$. For every constant $c>0$ and
$\epsilon>0$, there exists a  K\"ahler metric $\omega$ of
strictly positive constant scalar 
curvature on $X$ such that 
\begin{equation}
\label{eqfarein}
  \left | \frac {[\omega]\cdot F}{\sqrt{[\omega]^2}} - c \right |\leq \epsilon.
\end{equation}
\end{theointro}

This theorem can be thought of as an exploration of the K\"ahler cone
on certain del Pezzo surfaces. We are asking the question: which  
K\"ahler classes can be realized by K\"ahler metrics of strictly
positive constant scalar curvature?

\begin{remark} For the KE metric on $X_k$, $k=6, 7$, we have
the number $c$ takes the value
$$
c^\mathrm{ein}:=  \frac {[\omega^\mathrm{ein}]\cdot
  F}{\sqrt{[\omega^\mathrm{ein}]^2}} =  
\frac{c_1(X)\cdot F}{\sqrt{c_1^2(X)}} = \frac{2}{\sqrt{8-k}}.
$$
\end{remark}

\begin{remark}  It is an observation of LeBrun (see, for example,
\cite{LS}) that the hyperplane 
$$
H =\{\alpha\in H^2(X,\RR): c_1\cdot \alpha =0\}
$$
does not meet the K\"ahler cone of $X_k$ unless $k\geq 9$. It follows that if $k\leq 8$, every K\"ahler class $[\omega]$ satisfies $c_1(X)\cdot[\omega] >0$. Since this number represents, up to a factor of $4\pi$,  the integral of the scalar curvature over $X$, it follows that {\em any} CSC K\"ahler metric on $X_k$ must be positive, whatever the K\"ahler class.
\end{remark}

\vspace{12pt}
Now we turn to the case $k\geq 9$. In \cite{RS}, it was proved that there exist zero CSC K\"ahler metrics on $X_k$ in this case (the same result was previously known for $k\geq 14$ \cite{KLP}).  In the next result, we obtain CSC K\"ahler metrics of either sign:
\begin{theointro}
\label{fano}
There exists a $9$-point blow-up $X$ of $\CP^1\times \CP^1$ with no
non trivial holomorphic vector field and the
following properties: 
\begin{enumerate}
\item there is a scalar-flat K\"ahler metric on $X$;
\item for every $c\in \RR$ and $\epsilon>0$, there exists a metric of
constant scalar curvature $\omega$ on $X$ such that
\begin{equation}
\label{eqfarsfk}
  \left | \frac {[\omega]\cdot c_1(X)}{\sqrt{[\omega]^2}} - c \right |\leq \epsilon.
\end{equation}
\end{enumerate}
In addition, any further blow-up of $X$ admits  metrics with the same
properties.
\end{theointro}

Part (i) of this Theorem was proved in \cite{RS}. It plays
in this discussion an analogous  role to the one played by
K\"ahler-Einstein metrics in Theorem~\ref{cp9}, and corresponds to
$c=c^{\mathrm{sfk}}=0$. Again, the new content of this Theorem is that
there are metrics of constant scalar curvature that represent K\"ahler
classes arbitrarily far from the hyperplane $H$.
inside $H^2(X, \RR)$.

\begin{remark} K\"ahler metrics of constant scalar curvature admit a
good perturbation theory: if the class $[\omega]$ is represented by a CSC K\"ahler metric, then the same is true for all sufficiently close K\"ahler classes, at least if $X$ carries no non-zero holomorphic vector fields \cite{LSc}.  This applies to the KE metrics on $X_6$ and  $X_7$ to give non-Einstein CSC metrics in all classes sufficiently close to $c_1(X)$. It also applies to give CSC K\"ahler metrics of either sign on $X_9$, by perturbing a known zero CSC metric.  It seems very difficult, however, to extend this perturbation theory to give more global results like Theorems~\ref{cp9} or \ref{fano}.
\end{remark}

\subsection{Strategy} The Theorems~\ref{cp9} and \ref{fano} will be
deduced from an extension of the constructions invented in~\cite{RS}
and the recent gluing theorem of Arezzo-Pacard~\cite{AP} for CSC
K\"ahler metrics.  The zero-CSC K\"ahler metrics in \cite{RS} were
obtained by resolving the singularities of a zero-CSC orbifold. The
gluing theorem in that work applied only to produce zero-CSC K\"ahler
metrics. The results of \cite{AP} allow us to work with CSC orbifolds
of non-zero scalar curvature, which are, however, constructed in
almost exactly the same way as in \cite{RS}. Indeed, as will be
explained below, these orbifolds are twisted products of Riemann
surfaces, where now the scalar curvature of the factors can be chosen
arbitrarily.


\begin{remark}
In \S9(i) of their paper, Arezzo and Pacard make the following construction: take a copy of $X_k$ with no holomorphic vector field and Tian's K\"ahler--Einstein metric (so $k\leq 7$). Their gluing theorem gives a positive CSC K\"ahler metric on the blow-up of $X_k$, that is to say, on $X_{k+1}$.  If $k=7$, then this
proves Theorem~\ref{cp9} in the particular case $c=2$, $k=8$.

We should point out that in contrast to this approach, our construction does not use any result about the existence of K\"ahler-Einstein metrics on del Pezzo surfaces and may be considered to be
more elementary for this reason.
\end{remark}

\begin{remark}The gluing theorem of \cite{AP} is not specific to
complex dimension $2$.    Accordingly, in \S\ref{higher} we give some new but rather special higher-dimensional examples of negative CSC K\"ahler metrics.
\end{remark}


\subsection*{Acknowledgments} We thank Claudio Arezzo, Denis Auroux
and Frank Pacard for stimulating discussions. 

\section{The factory}
\label{secfactory}
We review the construction of scalar-flat K\"ahler orbifold metrics of~\cite{RS} and adapt
it to the case of K\"ahler  metrics of constant scalar curvature.  

\subsection{Orbifold Riemann surfaces}
\label{secgood}
We start with a closed Riemann surface $\overline 
\Sigma $ of genus $g$ with a 
finite set of orbifold points $P_1,P_2,\cdots, P_k$, with local
ramified cover of 
$q_1,q_2,\cdots, q_k >1$.
 Recall first the description of the fundamental group of the 
punctured Riemann surface $\Sigma=\overline\Sigma\setminus\{P_j\}$:
\begin{multline*}\label{e3.10.12.03}
\pi_1(\Sigma)= \\
\langle a_1,b_1,\ldots,a_g,b_g,l_1,\ldots, l_k:
[a_1,b_1][a_2,b_2]\ldots[a_g,b_g]l_1\ldots l_k = 1\rangle
\end{multline*}
Here the $a_j$ and $b_j$ are standard generators of $\pi_1(\overline\Sigma)$
 and  $l_j$ is
(the homotopy class of) a small loop around $P_j$. The orbifold
fundamental group is defined by
\begin{equation*}\label{e4.10.12.03}
\pi^\mathrm{orb}_1(\ovS) = \pi_1(\Sigma) / G
\end{equation*}
where $G$ is the normal subgroup  of $\pi_1(\Sigma)$ generated by  $l_1^{q_1}, l_2^{q_2},\cdots,
 l_k^{q_k}$.

The orbifold Euler characteristic is defined by
$$ 
\chi^\mathrm{orb}(\overline\Sigma):= \chi^\mathrm{top}(\ovS)  - \sum_{j=1}^k
(1- \frac 1 {q_j}) .  $$ The question of whether $\overline{\Sigma}$
carries a CSC orbifold K\"ahler metric was considered by Troyanov
\cite{Tr}. Let us call an orbifold Riemann surface ``good''
if its orbifold universal cover admits a compatible K\"ahler metric of
CSC $\kappa_1$, say.  
By the Gauss--Bonnet theorem (applied to the surface
with boundary obtained by removing a small disc around each of the
$P_j$), if $\ovS$ is good, then the sign of $\kappa_1$ is the same as
the sign of $\chi^{\mathrm{orb}}(\ovS)$. 
However, not every orbifold Riemann
surface is ``good'': the {\em tear-drop}, which is $S^2$ with one
orbifold point of any order $\geq 2$, is simply connected but admits no
compatible metric of constant curvature. 

The following summarizes the facts we shall need if the orbifold euler
characterstic is non-positive:

\begin{prop} The orbifold Riemann surface $\ovS$ is always good if
$$\chi^{\mathrm{orb}}(\ovS) \leq 0.$$
  Such $\ovS$ carries no
non-trivial holomorphic vector fields if $\chi^{\mathrm{orb}}(\ovS)<0$
or if $\chi^{\mathrm{orb}}(\ovS)=0$ and there is at least one orbifold
point.
\end{prop}
\begin{proof}The existence of metrics of constant scalar curvature is
contained in \cite{Tr}. The statement about holomorphic vector fields
follows either by lifting such a field to the universal cover or from
the orbifold version of the Hopf index theorem: a holomorphic vector
field must vanish at each of the orbifold points.
\end{proof}

For good orbifolds with  positive orbifold euler characteristic, the
universal cover must be $\CP^1$ and the orbifold fundamental group
must be one of the finite subgroups of $\SO(3)$:
\begin{prop} Let $\ovS$ be a good orbifold with
$\chi^{\mathrm{orb}}(\ovS)>0$. Then $\ovS$ is biholomorphic to
$\CP^1/\Gamma$, where $\Gamma$ is a finite subgroup of
$\SO(3)$. There are either $2$ marked points of the same order $q\geq
2$ or $3$ marked points with orders $\{2,2,q\}$, $\{2,3,3\}$,
$\{2,3,4\}$, $\{2,3,5\}$. The corresponding groups are cyclic of order
$q$, dihedral, tetrahedral, octahedral, or icosahedral respectively.

Furthermore, $\ovS$ carries no non-trivial holomorphic vector fields
if and only if $G$ is not cyclic.
\label{p1.17.12.4}
\end{prop}
\begin{proof} We need only discuss the last part. The Euler vector
field ($z\partial_z$ in affine coordinates) descends to the quotient
if $G$ is the cyclic group acting with fixed-points at $z=0$ and
$z=\infty$. On the other hand, no non-trival holomorphic vector field
on the sphere can have 3 zeros, so none of the other orbifolds can
admit such vector fields.
\end{proof}

\subsection{Desingularization of orbifold ruled surfaces}

From now we on assume that $\ovS$ is a good orbifold Riemann
surface and we endow $\ovS$ with an orbifold K\"ahler metric $\bar
g^\mathrm{\Sigma}$ of constant curvature $\kappa_1$.  Note that, just
as for ordinary Riemann surfaces, we have $$
\ovS = \cU/\pi_1^{\mathrm{orb}}(\ovS)
$$ 
where the fundamental group acts by isometries on the universal
cover $\cU$ which is equal to $\HH$, $\EE$ or $\CP^1$, according as
$\kappa_1$ is negative, zero, or positive.

Let $g^\mathrm{FS}$ be the Fubini-Study metric with curvature
$\kappa_2>0$, on $\CP^1$. 
Let $\rho:
\pi_1^\mathrm{orb}(\overline\Sigma) \rightarrow \SU(2)/\ZZ_2$ be
a homomorphism. Then $\pi_1^\mathrm{orb}(\ovS)$ acts
isometrically on $\cU\times\CP^1$. 
 We deduce a K\"ahler metric  $\bar g_\rho$ of constant scalar
curvature $s=2(\kappa_1+\kappa_2)$ on the orbifold quotient
$$\overline M_\rho = \ovS\times_\rho\CP^1 =
(\cU\times\CP^1)/ \pi_1^\mathrm{orb}(\overline\Sigma).$$

The idea now is to apply the results of \cite{AP} to obtain CSC
K\"ahler metrics on the minimal resolution of singularities,
$\widehat{M}_\rho$, say, of $\ovM_\rho$. We note that the constant $s$
can be chosen arbitrarily if $\ovS$ is hyperbolic, but must be
positive otherwise.  Here is the full statement, which is parallel to
Theorem~D of \cite{RS}.
\begin{theointro}
\label{maintheo}
Let $\ovS$ be a good compact orbifold Riemann surface as above,
carrying no non-trivial holomorphic vector fields. Suppose that
$\rho: \pi_1^\mathrm{orb}(\ovS ) \rightarrow \SU(2)/\ZZ_2$ is a
homomorphism that is {\em irreducible} in the sense that the induced action
of $\pi_1^\mathrm{orb}(\ovS)$ fixes no point of $\CP^1$. 
Equip $\ovM_\rho$ with a twisted product metric of CSC $s$ as above. 

Then the minimal resolution $\hM_\rho \to \ovM_\rho$ carries a CSC
K\"ahler metric with scalar curvature very close to $s$.

In addition, any further blow-up of $\widehat M_\rho$  carries a CSC K\"ahler
metric  with scalar curvature very close to $s$.
\end{theointro}

\begin{proof} We have already equipped $\ovM_\rho$ with a K\"ahler
orbifold metric of constant scalar curvature $s$. On the other hand,
for each finite cyclic subgroup $G$ of $U_2$ with the property that
$G$ acts freely on $S^3$, there is an asymptotically locally euclidean
zero-CSC K\"ahler metric on the minimal resolution \cite{CS, RS}. 
If $s=0$, we are in the situation of \cite{RS}, so we may as well
assume that $s\not=0$. We must check that the hypotheses of \cite{AP}
are satisfied, i.e.\ that $\ovM_\rho$ has no non-trivial holomorphic
vector fields.  But this follows as in \cite{RS}, for we have placed
ourselves in the situation where the base $\ovS$ has no non-zero holomorphic
vector fields.

Hence, we can indeed apply \cite{AP} to obtain CSC K\"ahler metrics
$g$ on $\hM_\rho$.  For sufficiently small choices of the gluing
parameter, $g$ is very close to $g_\rho$, and the scalar curvature of $g$
will be very close to the number $s$.

To get the last statement, we again apply \cite{AP}, this time gluing
a copy of the Burns metric at a smooth point.
\end{proof}

\begin{remark}
The term ``minimal'' applied here is really unfortunate. According to
its usual meaning in complex surface theory, it means that there are
no rational curves of self-intersection $-1$. In this case, however, 
the proper transform of any singular fibre of $\ovM_\rho$ is always a
$(-1)$-curve! The term ``minimal'' is justified by the fact that we
are taking the minimal resolution of each of the orbifold
singularities separately. Only after doing so is it possible to blow
down the $(-1)$-curves.
\end{remark}



\section{Stability of parabolic ruled surfaces}
\label{secparab}
In order to make the results of the previous section more useful, we need a way to generate representations 
$\rho: \pi_1^{\mathrm{orb}}(\ovS)  \longrightarrow \SU(2)/\ZZ_2$.
A way to do this was introduced in \cite{RS}, through the notion of a (stable) parabolic ruled surfaces. We recall the main definitions here:

A {\em geometrically ruled surface} $\chM$ is by definition a minimal complex
surface obtained as $\chM= \PP(E)$, where $E\rightarrow \widehat \Sigma
$ is a holomorphic vector bundle of rank $2$ over a \emph{smooth} Riemann surface 
$\widehat \Sigma$.
The induced map $\pi:\chM \rightarrow \widehat \Sigma$ is called the
{\em ruling}.

A parabolic structure on $\chM$ consists of the following data:
\begin{itemize}
\item A finite set of distinct points $P_1,P_2,\cdots,P_n$ in
  $\widehat \Sigma$;
\item for each $j$, a choice of point  $Q_j \in F_j = \pi^{-1}(P_j)$;

\item for each $j$, a choice of  weight $\alpha_j \in   (0,1)\cap\QQ$.
\end{itemize}
A geometrically ruled surface with a parabolic structure will be called a
{\em parabolic ruled surface}.

If $S\subset\chM$ is a holomorphic section of $\pi$, we define its
slope
$$\mu(S) = S^2 +\sum_{Q_j\not\in S}\alpha_j - \sum_{Q_j\in S}\alpha_j;
$$
we say that a parabolic ruled surface is \emph{stable} if for
every holomorphic section $S$, we have $\mu(S) >0$.

\subsection{Iterated blow-up of a parabolic ruled surface}
\label{secitbup}

Let $\chM$ be a parabolic ruled surface. We shall now define a
multiple blow-up $\Phi:\hM \to \chM$ which is canonically
determined by the parabolic structure of $\chM$.

In order to simplify the notation, suppose that the parabolic
structure on $\chM$ is reduced to a single point $P\in \widehat
\Sigma$; let $Q$ be the corresponding point in $F = \pi^{-1}(P)$ and let
$\alpha = \frac pq$ be the weight, where $p$ and $q$ are two coprime
integers, $0<p<q$. Denote the Hirzebruch--Jung continued fraction
expansion of $\alpha$ by
\begin{equation}\label{e1.844}
\frac pq =  \cfrac{1}{e_1-\cfrac{1}{e_2-\cdots\cfrac{1}{e_l}}};
\end{equation}
define also
\begin{equation}\label{e20.844}
\frac {q-p}{q} = \cfrac{1}{e'_1-\cfrac{1}{e'_2-\cdots
\cfrac{1}{e'_m}}}.
\end{equation}
These expansions are unique if, as we shall assume, the $e_j$ and $e'_j$ 
are all $\geq 2$.

We give here a construction the iterated blow-up $\hM$: the fiber $F$
has self-intersection $0$. The first step is to
blow up $Q$, to get a diagram of the form
\begin{equation}\label{e2.844}
\xymatrix{
{}\ar@{-}[rr]^{-1} && *+[o][F-]{}
\ar@{-}[rr]^{-1} &&{}
}
\end{equation}
By blowing up the intersection point of these two curves we get the
diagram
\begin{equation}\label{e3.844}
\xymatrix{
{}\ar@{-}[rr]^{-2} && *+[o][F-]{} 
\ar@{-}[rr]^{-1} && *+[o][F]{}
\ar@{-}[rr]^{-2} &&{} 
}
\end{equation}
Then we perform an iterated  blow-ups of one of the two intersection of the
only $-1$-curve in the diagram. Given $\alpha=p/q$, there is a unique way
(cf.~\cite[Proposition 2.1.1]{RS}) to choose at
each step which
point has to be blown-up  in order to get the following diagram
\begin{equation*}
\xymatrix{
{}\ar@{-}[r]^{-e_1} & *+[o][F-]{}
\ar@{-}[r]^{ -e_{2}} &  *+[o][F-]{}
\ar@{--}[r] &  *+[o][F-]{}
\ar@{-}[r]^{-e_{l-1}} &  *+[o][F-]{}
\ar@{-}[r]^{-e_l} &  *+[o][F-]{}
\ar@{-}[r]^{-1} &  *+[o][F-]{}
\ar@{-}[r]^{-e'_{m}} &  *+[o][F-]{}
\ar@{-}[r]^{-e'_{m-1}} &  *+[o][F-]{}
\ar@{--}[r] &  *+[o][F-]{}
\ar@{-}[r]^{-e'_{2}} &  *+[o][F-]{}
\ar@{-}[r]^{ -e'_{1}} & 
}
\end{equation*}
where the $-e_1$-curve is the proper transform of the original fiber
$F$.

More generally, if $\chM$ has more parabolic points, we perform the same operation
for every point and get a corresponding iterated blow-up
$\Phi:\hM\to\chM$. 

\subsection{The theorem of Mehta-Seshadri}
We now give a practical way of using Theorem~\ref{maintheo}. Given a
parabolic ruled surface $\chM \to \hS$,
 we deduce a orbifold Riemann surface $\ovS$ by introducing a orbifold
 singularity of order $q_j$ at every parabolic point $P_j\in \hS$ of weight
 $p_j/q_j$. As a corollary of Mehta-Seshadri theorem~\cite{MS}, we have the
 following proposition.
\begin{prop}
Let $\chM\to\hS$ be a parabolically stable ruled surface. Then
there exists an irreducible representation $\rho :
\pi_1^\mathrm{orb}(\ovS)\to \SU(2)/\ZZ_2$ such that $\hM \simeq
\hM_\rho$, where $\hM$ is the iterated blow-up of the parabolic ruled
surface $\chM$ as described in Section~\ref{secitbup} and $\hM_\rho$ is the
smooth resolution of the orbifold $\ovM_\rho$ defined in Section~\ref{secfactory}.
\end{prop}
\begin{proof}
This is a direct consequence of \cite[Theorem 3.3.1]{RS}.
\end{proof}

In conclusion we can reformulate Theorem~\ref{maintheo} as follows.

\begin{theointro}
\label{maintheo2}
Let $\chM\to \hS $ be a parabolically stable ruled surface, with parabolic
weights $\alpha_j = p_j/q_j$ where $p_j$ and $q_j$ are positive coprime integers.
Let $\ovS$ be the orbifold ruled surface deduced from $\hS$ according
to the parabolic structure.

Then the conclusions $(i)$ and $(ii)$ of Theorem~\ref{maintheo} hold
with $\hM_\rho$ replaced by $\hM$.
\end{theointro}

\section{The goods: proof of Theorems~\ref{cp9} and \ref{fano}}
\subsection{A stable parabolic bundle}\label{secaparab}
Let $\ovS$ be $\CP^1$ with 
$3$ orbifold points $P_1$, $P_2$, $P_3$, each of weight $2$. Then the
orbifold Euler characterstic is $1/2$ and $\ovS$ is good, and we can
equip $\ovS$ with a metric of CSC $\kappa_1>0$.  
Now, we consider the ruled surface $\pi:\CP^1\times\CP^1\to\CP^1$
where $\pi$ is, say, the projection on the first factor.  We pick a
point $Q_j$  in the fibre $\pi^{-1}(P_j)$, and give it weight
$\alpha_j=1/2$.

\begin{lemma}
  For a generic choice of points $Q_1,Q_2,Q_3$, the  parabolic
  ruled surface $\chM\to\CP^1$ defined above is parabolically stable.
\end{lemma}
\begin{proof}
We just need to arrange so that
 $2$  points $Q_j$ never belong
to the same constant section of $\chM\to \CP^1$. 
\end{proof}

\subsection{Conclusion} The iterated blow-up $\hM$ of Theorem~D in
this case consists of performing 2 blow-ups in each fibre, so that
$\hM$ is a $6$-point blow-up of $\CP^1\times \CP^1$. Thus this result
yields positive-CSC K\"ahler metrics on
$X_6$.  By adjusting $\kappa_1$ and $\kappa_2$ and taking the gluing
parameter to be small enough, we can get any positive value of $c$ in
\eqref{eqfarein}.

To deal with $X_7$ and $X_8$, we can use the last part of
Theorem~\ref{maintheo} to perform further blow-ups on the previous
example. Alternatively, we can replace one of the parabolic weights by
$1/3$ or $1/4$ in the above.  The complex surface $X_8$ can also be
constructed by adding a fourth parabolic point $P_4$, and taking all
weights equal to $1/2$. It is readily seen that for generic choices of
the points $Q_j$, the resulting parabolic ruled surfaces are stable.

\begin{remark} We deliberately chose $\kappa_1$ and $\kappa_2$ so as
to get CSC metrics that are far from being Einstein. However, if we
take $\kappa_1=\kappa_2$ then it seems possible that with a little
more work this gluing construction might be adapted to give the KE
metrics for (certain complex structures on) $X_6$ and $X_7$! From this
point of view, it would also be interesting to try to extend the
analysis to allow $\ovM_\rho$ to carry non-trival holomorphic vector
fields. As the analysis in \cite{LSc} shows, it is quite possible for
the deformation problem to be unobstructed in this case, although the
whole story is then much more delicate.
\end{remark}

\subsection{Proof of Theorem~\ref{fano}} The proof follows the same
lines as in \cite{RS} (or as above), taking $\ovS$ to be $\CP^1$ with
4 orbifold points of order $2,2,2,3$. The orbifold euler
characteristic is $-1/6$ and so we can equip $\ovM_\rho$ with a CSC
K\"ahler metric with normalized scalar curvature
$c_1\cdot[\omega]/\sqrt{[\omega]^2}$ equal to any given real
number. Applying Theorem~D, we get the result claimed in
Theorem~\ref{fano}. 

\section{Explicit representations}

Our use of the theory of stability to generate representations of the
orbifold fundamental group results in a certain loss of
explicitness---one of the merits of gluing constructions is that they
give ``nearly'' explicit metrics, in the sense that for small values
of the gluing parameter, the glued metric is close to the original
one.  From this point of view, it seems desirable to try to make
explicit these representations, at least in a few examples, and this
section is devoted to such a task.

\subsection{Platonic orbifold surfaces}

Let $\Gamma$ be one of the finite sub\-groups of $\SO(3)$ as in
Proposition~\ref{p1.17.12.4}. Then there is a canonical orbifold
$\ovM=\CP^1\times\CP^1/\Gamma$, where $\Gamma$ acts diagonally on the
product. If $\Gamma$ is not cyclic, then this orbifold will not
support any non-trivial holomorphic vector fields and the resolution
will carry a positive CSC K\"ahler metric.

\subsection{Representations of the fundamental group of the $\mathbf{4}$-fold
punctured sphere}
\label{secexamples}
Let $\ovS$ be the orbifold Riemann surface used in the proof of
Theorem~\ref{fano}. We shall construct representations of the orbifold
fundamental group of $\ovS$.  Let the points be denoted $P_j$, with
orders $2$, $2$, $2$, $3$ respectively, and let
$l_j$ denote the homotopy class of a small loop around $P_j$.
The orbifold fundamental group is
$$
\pi_1^\mathrm{orb}(\ovS) = \langle l_1,l_2,l_3,l_4 \,:\, l_1l_2l_3l_4 = 
l_1^2 = l_2^2 = l_3^2 = l_4^3 =  1 \rangle.
$$
Identifying $\SU(2)$ with the group of unit quaternions, we write down a
family of 
representations $\pi_1^\mathrm{orb} \to \SU(2)/\ZZ_2$
$$
\rho(l_1) = \pm i, \rho(l_2) = \pm j, \rho(l_3)= \pm k
\exp(\pi(i\cos \phi + j\sin \phi)/3),\hspace{2cm}$$
$$
\hspace{2cm} \rho(l_4) = \pm \exp(-\pi(i\cos \phi + j\sin \phi)/3)
$$
Denote by $R_1$, $R_2$ and $R_3$ $R_4$ the corresponding rotations of
$\RR^3$. Then $R_1$, $R_2$ and $R_3$, while $R_4$ is a
$2\pi/3$-rotation,  see Figure~\ref{figrot}.

\ifpix
\begin{figure}[ht]
\epsfig{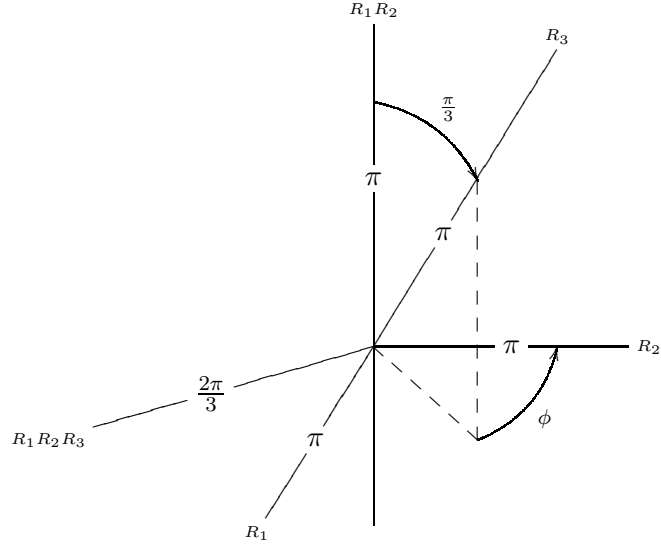}\medskip
\caption{A representation of the fundamental group of the 
 $4$-fold punctured sphere \label{figrot}}
\end{figure}
\fi

This representation is irreducible and so gives rise to a CSC K\"ahler orbifold
$\ovM_{\rho}$. The corresponding minimal resolution $\hM_{\rho}$ is
the 9-point blow-up of a minimal ruled surface $\chM \to \CP^1$. We do
not know which ruled surface this is for general $\phi$, but if
$\phi=0$, we 
sketch an argument explaining why $\chM$ is isomorphic to $\CP^1\times
\CP^1$, so that the scalar-flat K\"ahler manifold $\hM_\rho$ is
isomorphic to a $9$-point blow-up of $\CP^1\times \CP^1$.  

Taking $\phi=0$ from now on, we note that
the $x$-axis is stable under the
action of $\rho$: it is fixed point-wise by $R_1$ and $R_4$ and 
flipped by $R_2$ and $R_3$. Set $Q_1 = (1,0,0)$, $Q_2 = (-1,0,0)$ and
consider the Riemann surface
$$
S := (\CP^1\setminus \{P_j\}) \times_\rho \{Q_1,Q_2\}.
$$
This sits in $\ovM$, and its 
closure $\overline S$ is
a $2$-fold ramified cover $p:\overline
S\rightarrow \ovS$. In $\overline{S}$, there are two
orbifold singularities sitting over each of $P_1$ and $P_4$; on the
other hand the cover is ramified over $P_2$ and $P_3$ and the
corresponding points in $\overline{S}$ are smooth. Thus $\overline S$
is an orbifold version of $\CP^1$, with two singular points of weight
2 and two singular points of weight 3.
Denote by $\overline{N}$ the associated orbifold ruled surface, and by
$\widehat{N}$ its minimal resolution.

Downstairs in $\overline{M}$ there are $8$ orbifold singularities,
say,
$$
A_j, B_j \in \ovpi^{-1}(P_j).
$$
Upstairs in $\overline{N}$, $P_1$ and $P_4$ have been double-covered
by points $P_1^{\pm}$, $P_4^{\pm}$ and there are accordingly
singularities
$$
A_1^{\pm}, B_1^{\pm} \in \ovpi^{-1}(P_1^{\pm}),\;\;
A_4^{\pm}, B_4^{\pm} \in \ovpi^{-1}(P_4^{\pm}).
$$
The double cover is ramified over $P_2$ and $P_3$, so the points 
$A_2'$, $A_3'$, $B_2'$, $B_3'$
corresponding to $A_2$, $A_3$, $B_2$, $B_3$ are smooth. There
is a commutative diagram
$$
\xymatrix{
 \overline N \ar[d]_{p^*\bar\pi}\ar[dr] & & \widehat N\ar@{-->}[d]
 \ar[ll] \\
  *++{\overline S} \ar@{^(->}[r] \ar[dr]_p & \ar[d]^{\bar\pi}
  \overline M & \ar[l] \widehat M \\
& \ovS  &
}
$$
where the dashed arrow is a meromorphic map which becomes a smooth
ramified double cover if the points $A_2'$, $A_3'$, $B_2'$, $B_3'$ are
blown up.

The reason for passing to $\overline N$ is that it contains two
parallel
sections $s_1$ and $s_2$, defined by $s_1(Q_1)=Q_1$, $s_1(Q_2)=Q_2$ and
$s_2(Q_1)=Q_2$, $s_2(Q_2)=Q_1$; let ${\overline S}_j$ denote the closure
of $s_j(\overline S)$ in $\overline N$. Because the $s_j$ are
parallel, ${\overline S}_j^2=0$. By a computation that we omit,
the 
proper transforms $\widehat S_j$ of $\overline S_j$  in  $\widehat N$
verify $\widehat S_1\cdot \widehat S_2 = 0$ and  
$$ \quad \widehat S_1^2 = \widehat
S^2_2 =  -2,
$$
see Figure~\ref{fig2cover}.

\ifpix
\begin{figure}[ht]
\epsfig{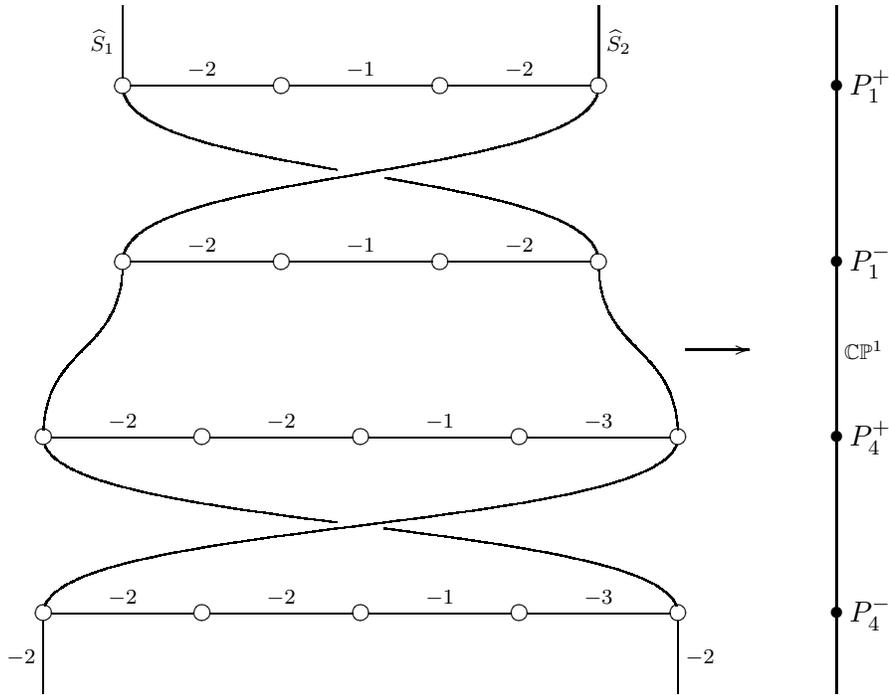}\medskip

\caption{Configuration of curves in $\widehat N$\label{fig2cover}}
\end{figure}
\fi

\noindent 
After a sequence of $5$ blow-downs in a suitable order in $\widehat
M$ (corresponding to  $10$ blow-downs on $\widehat N$) we obtain the
configuration of curves shown in 
Figure~\ref{figbdown}:
\ifpix

\begin{figure}[ht]
\epsfig{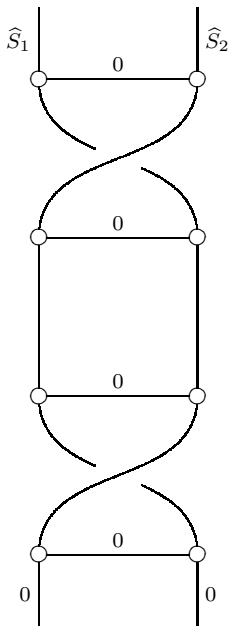}\medskip

\caption{10-point blow-down of $\widehat N$ \label{figbdown}}
\end{figure}
\fi
It is a minimal ruled 
surface with a divisor $D$ verifying $D\cdot F=1$ and $D^2=0$, where
$F$ is the class of a generic fibre.
Therefore it is isomorphic to  $\CP^1\times\CP^1$. Thus
$\widehat M$ is a five-point blow-up of the minimal resolution 
of $(\CP^1\times\CP^1)/\ZZ_2$, where 
the action of $\ZZ_2$ rotates each factor by an angle $\pi$.

It is easy to see that the latter complex surface is a four-point
iterated blow-up of $\CP^1\times\CP^1$.
Consider the proper transforms $\widehat S_1$ of
$\overline S_1= (0\times
\CP^1)/\ZZ_2$ and  $\widehat S_2$ of
$\overline S_2= (\CP^1 \times 0)/\ZZ_2$. Then we have 
$\widehat S_1\cdot \widehat S_2=0$ and $\widehat S_1^2 = \widehat
S_2^2 =  -1$. 
Hence we have the configuration of curves shown in Figure~\ref{fig2cov2}
and after $4$ blow-downs in a suitable order, we obtain
  $\CP^1\times\CP^1$.

\ifpix
\begin{figure}[ht]
\epsfig{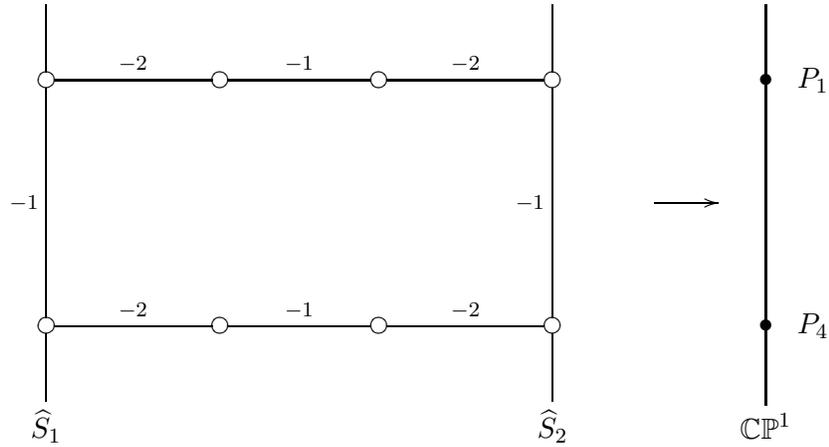}\medskip

\caption{$\widehat M$ with five $(-1)$-curves blown down  \label{fig2cov2}}
\end{figure}
\fi

\subsection{Representations of the fundamental group of the 
punctured torus}

Let $\TT$ be a compact elliptic curve with one marked point $P$ of
weight $2$, let $\ovS$ stand for the corresponding orbifold Riemann
surface. In terms of standard generators,
$$
\pi_1^\mathrm{orb}(\ovS) = \langle a,b,l\,:\, [a,b]l = l^2 = 1\rangle.
$$
We define a representation $\rho: \pi_1^\mathrm{orb}(\ovS) \to \SU(2)/\ZZ_2$ 
as follows
$$ \rho(a)= \pm e^{\pi k/4}i,\; \rho(b) = \pm j, \; \rho(l)= \pm k.
$$
The corresponding rotations $R_a$, $R_b$, $R_l$ are shown in
in Figure~\ref{figrot2}.

\ifpix
\begin{figure}[ht]

\epsfig{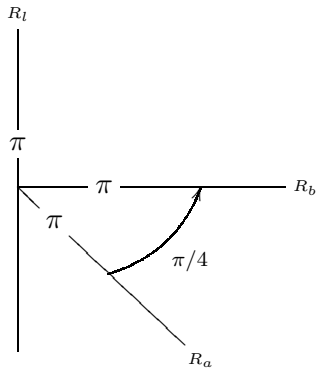}\medskip

\caption{A representation of the punctured torus \label{figrot2}}
\end{figure}
\fi

Denote by $\ovM$ the CSC K\"ahler orbifold ruled surface arising from this
representation, by $\hM$ the minimal resolution, 
 and by $\chM\to\TT$
a  minimal model for $\hM$.

Then we have
\begin{prop}
\label{propt}
The surface $\hM$ is a double blow-up of $\chM$, and carries a CSC K\"ahler metric.
On the other hand, no smooth blow-down of
$\widehat M$ admits CSC K\"ahler metric.
\end{prop}

We give a sketch of the proof of this.
The existence of a metric on $\widehat M$ is assured by
Theorem~\ref{maintheo}. Consider, then, the non-existence statement.
The geometrically ruled surface $\chM$ cannot admit a CSC K\"ahler metric for
$c_1^2=0$, but $c_1\not=0$. In $\hM$ the only curve that can be blown
down smoothly is in the fibre over $P$. Denote by $X$ this blow-down.

As in the previous example, we can pass to a double cover by
remembering the two points on the unit sphere on the axis of $R_l$. We
obtain a double-cover $\overline N$ of $\ovM$ and a corresponding minimal
resolution $\widehat N$. By considering the geometry of this double-cover,
one can show that the double-cover $X'$ of $X$ is a 2-point blow-up of
$\PP(L_1\oplus L_2) \to \TT$, where the points now lie on the sections
corresponding to $L_1$ and $L_2$. Such $X'$ does not admit a CSC K\"ahler
metric by \cite[Prop.~3.1]{LS}, and so $X$ cannot admit a CSC K\"ahler metric
either.

\section{Higher-dimensional examples}\label{higher}

Since the results of \cite{AP} apply in all dimensions, it is natural to try to extend the foregoing results to higher dimensions.  The main problem here is that very little appears to be known about the following basic
\begin{question}  If $G$ is a finite subgroup of $U_{m}$ acting freely on $\CC^{m}\setminus 0$, and $X\to \CC^{m}/G$ is a resolution of singularities, does there exist an asymptotically locally euclidean zero CSC K\"ahler metric on $X$?
\end{question}

Joyce \cite{J} proved that if $X$ can be chosen to be a crepant resolution ($c_{1}(X)=0$), then $X$ does carry an ALE Ricci-flat metric.  In dimensions $2$ and $3$, this is the case if and only if $G$ (is conjugate to) a subgroup of 
$\SU(m)$.  In higher dimensions, however, it is not known in general which singularities $\CC^{m}/G$ admit crepant resolutions.

The other class of examples come from explicit constructions of K\"ahler metrics on total spaces of line-bundles, as described, for example, in \cite{HS} which contains an extensive survey of the literature.  The following result follows easily from the methods of that paper, and is probably well known to many.  We sketch the proof in order to make this paper self-contained.

\begin{theo}  \label{thmHS}
If $G = \{1,\omega, \cdots, \omega^{k-1}\}$, where 
$\omega = e^{2\pi i/k}$ then there is a resolution $X$ of $\CC^{m+1}/G$ that admits an asymptotically locally euclidean zero-CSC  K\"ahler metric.
\end{theo}
\begin{proof} We note that for this group $G$,  we can take $X$ to be the total space of $\cO(-k) \to \CP^{m}$. Indeed, the total space of $\cO(-1)$ can be represented as  $U/\CC^{*}$, where
$$
U = \{(w,z_{0},\ldots, z_{m})\in \CC\times (\CC^{m}\setminus 0)\}$$
the action of $\CC^*$ is given by
$$t\cdot(w,z_{0},\ldots, z_{m}) =
(t^{-1}w,t z_{0},\ldots, t z_{m}), 
$$
and the map
$$
\beta(w,z_{0},\ldots, z_{m}) \mapsto (wz_{0}, wz_{1}, \ldots, wz_{m})
$$
is an invariant description of the blow-up of the origin of $\CC^{m+1}$.   
The action of the group $G$ lifts to the action
$$
\omega^{r}(w,z_{0}, z_{1},\ldots, z_{m}) = 
(\omega^{r}w, z_{0}, z_{1},\ldots, z_{m})
$$
on $U$, with quotient equal to the total space of $\cO(-k)$. 

Following an Ansatz that goes back to Calabi, we
seek  a K\"ahler metric on $X$ of the form
\begin{equation}\label{e2.8.4.5}
\omega_f = \omega_0  + i\del\delb f(t)
\end{equation}
where $\omega_{0}$ is the standard K\"ahler form on $\CP^{m}$ and $t$ is logarithm of the fibre-distance function.  

To be more explicit, we work with affine coordinates $(z_{1}, \ldots, z_{m})$ on $\CP^{m}$ so that
$$
\omega_0 = (i/2)\del\delb\log(1 + |z_1|^2 + \cdots + |z_m|^2).
$$
In the same coordinates, the fibre metric on $\cO(-k)$ is just given by
$$
h^{k} = (1+ |z_1|^2 +\cdots + |z_m|^2)^k.
$$
Denote by $w$ a holomorphic fibre coordinate, set $z = \log w$ so that
\begin{equation}\label{e1.8.4.5}
2t = z + \zb + k \log h
\end{equation}
We have
$$
i\del\delb t = (ik/2)\del\delb \log h = k\omega_0
$$
so that \eqref{e2.8.4.5} can be written
$$
\omega_{f}=(1 + kf'(t)) \omega_{0} + f''(t)i\del t\wedge \delb t.
$$

As explained in \cite{HS}, it is better to pass to momentum coordinates, so we introduce
\begin{equation}\label{e3.8.4.5}
\tau = f'(t)\mbox{ and }\varphi(\tau) = f''(t).
\end{equation}
The advantage of this is that the scalar curvature is given by a very simple formula which we now derive.  First of all,
\begin{equation}\label{e1.9.4.5}
\omega_f^{m+1} = (1+k\tau)^{m}h^{-m-1}\varphi(\tau)\psi\wedge\overline{\psi}
\end{equation}
where $\psi$ is a holomorphic $(m+1)$-form, so that
\begin{align}\label{e2.9.4.5}
\rho(\omega_f) & = -i\del\delb \log[(1+k\tau)^{m}h^{-m-1}\varphi(\tau)]\\
& =2(m+1)\omega_0 -i\del\delb \log[(1+k\tau)^{m}\varphi(\tau)] \nonumber.
\end{align}
Now for any function of $u(\tau)$,
\begin{equation}\label{e2.5.9.4.5}
i\del\delb u(\tau) =  k\varphi u'\omega_0 + \varphi(\varphi u')'i\del t\wedge \delb t
\end{equation}
where prime denotes differentiation with respect to $\tau$. Applying this with
$$
u(\tau) = m\log(1+k\tau) +\log \varphi(\tau)
$$
we obtain
\begin{equation}\label{e3.9.4.5}
\rho(\omega_f) = [2(m+1) - k\varphi u']\omega_0 - \varphi(\varphi u')' i\del t\wedge \delb t.
\end{equation}
The scalar curvature $\sigma(\omega_f)$ is given (up to a factor of $2$) by
$$
\sigma(\omega_f)\omega_f^{m+1} = (m+1)\rho(\omega_f)\wedge \omega_f^m
$$
which leads to the final formula
\begin{equation}\label{e4.8.4.5}
\sigma(\omega_f) = \frac{2m(m+1)}{1+k\tau} -
\frac{1}{(1+k\tau)^m}\frac{\rd^2}{\rd \tau^2}\left((1+k\tau)^m\varphi(\tau)\right).
\end{equation}

Thus $\sigma(\omega_f)=0$ if
\begin{equation}\label{e2.7.4.5}
(1+k\tau)^m \varphi(\tau) = \frac{2}{k^2}(1+k\tau)^{m+1} + a\tau +b
\end{equation}
where $a$ and $b$ are constants of integration.  The boundary conditions for smoothness of the metric at the zero-section are 
\begin{equation}\label{e3.7.4.5}
\varphi(0) = 0, \varphi'(0) = 2
\end{equation}
which gives
\begin{align*}
P(\tau) :&= (1+k\tau)^m\varphi(\tau) \\
&= \frac{2}{k^2}\left(
(1+k\tau)^{m+1} +(k-m-1)(1 +k\tau) + m - k\right).
\end{align*}
Geometrically, we have repackaged the complex line-bundle $X$ as the total space of an $S^{1}$-bundle over $\CP^{m}\times[0,\infty)_{\tau}$; to check that $\omega_{f}$ is really a metric on this $S^{1}$-bundle we need that $\varphi(\tau)>0$
for all $\tau\in(0,\infty)$. To see that this is the case, note that
$$
P'(\tau) = 0 \mbox{ if and only if } (1 + k\tau)^{m} = 1 -\frac{k}{m+1}
$$
and this latter equation cannot be satisfied if $\tau>0$. Since $P'(0)=2$, we see that $P'(\tau)>0$ for all $\tau>0$  and since $P(0)=0$ we have $P(\tau)>0$ for all $\tau>0$ as claimed.  
Thus our zero-CSC K\"ahler metric on $X$ corresponds to the function
\begin{multline}\label{e6.9.4.5}
\varphi(\tau) =
\frac{2}{k^2}\Big(
(1+k\tau) + 
(k-m-1)(1 +k\tau)^{1-m} \\ + (m - k)(1+k\tau)^{-m} \Big).
\end{multline}

Finally we check the asymptotics of this metric.  For this, we return to
the blow-up map at the beginning of the proof.  In the present coordinates, this is
\begin{equation}
x_0 = w^{1/k}, x_1 = w^{1/k}z_1,\ldots x_m = w^{1/k}z_m,
\end{equation}
where $(x_{0}, x_{1}, \ldots, x_{m})$ are standard linear coordinates on $\CC^{m+1}$, and the fractional power of $w$ corresponds to the passage from $\cO(-k)$ to its $k$-fold cover $\cO(-1)$.
The standard K\"ahler form on $\CC^{m+1}$ is
\begin{equation}
\eta = i\del\delb (|x_0|^2 + \cdots + |x_m|^2) = i\del\delb (|w|^{2/k}h) = i\del\delb ae^{2t/k}.
\end{equation}
Thus, this standard K\"ahler form is of the form \eqref{e2.8.4.5} with
\begin{equation}
f(t) = f_0(t) = e^{2t/k} - t/k.
\end{equation}
This yields
\begin{equation}\label{e5.9.4.5}
1+k \tau_0 = 2e^{2t/k},\;\; \varphi_0(\tau_0) = \frac{2}{k^2}(1+k\tau_{0}).
\end{equation}

Now since
\begin{equation}\label{e7.9.4.5}
dt = \frac{d\tau_0}{\varphi_0(\tau_0)} = \frac{d\tau}{\varphi(\tau)}
\end{equation}
we see that by suitable choice of integration constants, $\tau-\tau_0$ is order $\tau_0^{1-m}$ for large $\tau_{0}$, and so
$$
f(t) = f_0(t) +O(e^{2(1-m)t/k}) = f_0(t) + O( |x|^{2(1-m)})\mbox{ for large }|x|.
$$
Hence for large $|x|$,
\begin{equation}
\omega_{f} = i\del\delb[ (|x_0|^2 + \cdots + |x_m|^2) + O (|x|^{2(1-m)})],
\end{equation}
verifying that $\omega_{f}$ is ALE as required.
\end{proof}

We give two applications of this result:

\begin{theo} Let $S_{1}, \ldots, S_{m}$ be compact Riemann surfaces, each admitting a holomorphic involution $\iota_{j}$ with isolated fixed points. Suppose further that none of the $S_{j}$ is $\CP^{1}$, and that at least one is hyperbolic.  Let $\iota = (\iota_{1}, \ldots, \iota_{m})$  be the product involution on $S =S_{1}\times S_{2}\times \ldots \times S_{m}$.   Then the minimal resolution $\widehat{M}$ of 
$S/\iota$ admits a negative CSC K\"ahler metric.

In addition, any further blow-up of $\widehat M$  carries a CSC K\"ahler
metric  with negative scalar curvature.
\end{theo}

\begin{proof} Equip $S_{j}$ with a metric of constant curvature $\kappa_{j}$. The product metric on $S$ is not Ricci-flat because there is at least one hyperbolic factor.  Moreover the orbifold $S/\iota$ has no non-zero holomorphic vector fields. 
This is clear if each factor is hyperbolic, for then $S$ has no
non-zero holomorphic fields.   If one or more factor in the product is
flat, then $M$ does have non-zero holomorphic vector fields, but none
of them is preserved by $\iota$.   Hence for any choice
$(\kappa_{1},\ldots, \kappa_{m})$,  $S/\iota$ is an unobstructed CSC
K\"ahler orbifold, and every singularity is modelled on $\CC^{m}/(\pm
1)$.  Theorem~\ref{thmHS}, with $k=2$, now gives a resolution of this
singularity, and applying the desingularization theorem of \cite{AP}
gives the result.  The scalar curvature of $\widehat{M}$ is negative
since the sum of the $\kappa_{j}$'s is negative. 

For the last statement, the CSC K\"ahler metric on the blow-up is
obtained by gluing 
the ALE metric provided by Theorem~\ref{thmHS}, with $k=1$.
\end{proof}

\begin{theo}
Let $\ovS$ be a good compact orbifold Riemann surface as explained in Section~\ref{secgood}, carrying no non-trivial holomorphic vector fields, and with orbifold points of order $2$ only. Suppose that 
$\rho: \pi_1^\mathrm{orb}(\ovS ) \rightarrow \SU(2)/\ZZ_2$ is a 
homomorphism that is {\em irreducible} in the sense that the induced action
of $\pi_1^\mathrm{orb}(\ovS)$ fixes no point of $\CP^1$.  

Using the diagonal action of $\pi_1(\ovS)$ on $(\CP^1)^m$ induced by
$\rho$, define the twisted product $\ovM_\rho =\ovS\times_\rho
(\CP^1)^m$. This is an complex orbifold with singularities modelled on
$\CC^{m+1}/\{\pm 1\}$. 
Equip $\ovM_\rho$ with a twisted product metric of CSC. 
Then the minimal resolution $\hM_\rho \to \ovM_\rho$ carries a CSC
K\"ahler metric.

In addition, any further blow-up of $\widehat M_\rho$  carries a CSC K\"ahler
metric.  
\end{theo}
\begin{proof}
Again we will use the ALE spaces of Theorem~\ref{thmHS}, with $k=2$
 and $k=1$ (for the last statement).
 In order to use the gluing theorem of \cite{AP}, the only thing
 that needs to be checked is that $\ovM_\rho$ does not 
 carry any non-trivial holomorphic vector field. But, the argument is a
 straightforward generalization of the one given \cite{RS} in the case $m=1$.
\end{proof}

\begin{remark} In contrast to the previous theorem, by 
adjusting the curvatures of the base and the fibre here, in the case where $\chi^\mathrm{orb}(\ovS)<0$,  
we can arrange for the scalar curvature of $\widehat{M}_{\rho}$
to be positive, zero, or negative.
\end{remark}

\end{document}